\documentclass[final,10pt]{amsart}
\usepackage{amsmath,amsfonts,amssymb,amscd}
\usepackage[margin=0.4in]{geometry}
\usepackage{fullpage}
\usepackage{enumerate}

\usepackage{bm}
\usepackage{graphicx}
\usepackage{enumitem} 
\usepackage[all]{xy}
\usepackage{tikz}
\usepackage{tikz-cd}
\usetikzlibrary{arrows,tikzmark,calc}

\numberwithin{equation}{section}

\theoremstyle{definition}
\newtheorem{definition}[equation]{Definition}

\DeclareMathOperator\Aut{Aut}

\DeclareMathOperator\coh{coh}
\DeclareMathOperator\chr{char}

\DeclareMathOperator\End{End}

\DeclareMathOperator\gldim{gldim}
\DeclareMathOperator\GKdim{GKdim}
\DeclareMathOperator\grmod{grmod}
\DeclareMathOperator\GrMod{GrMod}
\DeclareMathOperator\gr{gr}
\DeclareMathOperator\id{id}

\DeclareMathOperator\Maxspec{Maxspec}

\DeclareMathOperator\Qcoh{Qcoh}
\DeclareMathOperator\qgrmod{qgrmod}
\DeclareMathOperator\QGrMod{QGrMod}

\DeclareMathOperator\supp{supp}
\DeclareMathOperator\Span{Span}
\DeclareMathOperator\Spec{Spec}
\DeclareMathOperator\tors{tors}
\DeclareMathOperator\Tors{Tors}

\DeclareMathOperator\trdeg{trdeg}
\DeclareMathOperator\trace{trace}

\newcommand\NN{\mathbb N}

\newcommand\ZZ{\mathbb Z}

\newcommand\cA{\mathcal A}

\newcommand\cC{\mathcal C}

\newcommand\cJ{\mathcal J}

\newcommand\cO{\mathcal O}
\newcommand\cR{\mathcal R}
\newcommand\cS{\mathcal S}
\newcommand\cX{\mathcal X}

\newcommand\fm{\mathfrak m}
\newcommand\fsp{\mathfrak p}

\newcommand\fsl{\mathfrak{sl}}

\newcommand\ba{\mathbf a}
\newcommand\bb{\mathbf b}

\newcommand\be{\mathbf e}

\newcommand\bp{\mathbf p}
\newcommand\bq{\mathbf q}

\newcommand\bsigma{\bm{\sigma}}
\newcommand\btau{\bm{\tau}}
\newcommand\bzero{\mathbf 0}

\newcommand\kk{\Bbbk}

\newcommand\fwa{W_1(\kk)}

\newcommand\wanr{W_n(R)}
\newcommand\wan{W_n(\kk)}
\newcommand\qwa{W_1^q(\kk)}
\newcommand\pwa{W_1^p(\kk)}

\newcommand\cnt{\mathcal Z}
\newcommand\inv{^{-1}}
\newcommand\iso{\cong}

\newcommand\tensor{\otimes}

\newcommand\grp[1]{\langle #1 \rangle}

\renewcommand{\to}{\ensuremath{\longrightarrow}}

\title{The Weyl algebra and its friends: a survey}
\author[Gaddis]{Jason Gaddis}
\address{Miami University, Department of Mathematics, Oxford, Ohio 45056} 
\email{gaddisj@miamioh.edu}


\begin{document}

\begin{abstract}
We survey several generalizations of the Weyl algebra including generalized Weyl algebras, twisted generalized Weyl algebras, quantized Weyl algebras, and Bell-Rogalski algebras. Attention is paid to ring-theoretic properties, representation theory, and invariant theory.
\end{abstract}

\maketitle

\section{Introduction}

The Heisenberg Lie algebra $H$ is generated over a field by $x,y,z$ with Lie bracket $[x,y]=z$ and $[z,x]=[z,y]=0$. In the enveloping algebra $U(H)$ the image of $z$ is central and hence $(z-1)$ is a two-sided ideal. The quotient $U(H)/(z-1)$ is now called the \emph{Weyl algebra}, a termed originally coined by Dixmier \cite{dix}. But the algebraic theory of this object goes back as far as to Born, Dirac, Heisenberg, Hirsch, and Littlewood, among others. For a historical account, the reader is directed to the survey of Coutinho \cite{Cout}. Volumes have been written on the Weyl algebra along with its properties and representation theory. In this survey we will give only passing acknowledgement to these and our focus will be on its generalizations. 

Let $\kk$ be an algebraically closed field of characteristic zero\footnote{These hypotheses are necessary for parts of the discussion but not all. The reader is encouraged to follow up with the given references (of which there are many) for minimal hypotheses on the field for particular results.}. Now and forevermore, rings (and algebras) are assumed to be associative. A module over a (noncommutative) ring $R$ is a left module. If $R$ is a ring and $a,b \in R$, then the commutator is denoted $[a,b]=ab-ba$. We denote the center of $R$ by $\cnt(R)$. For a positive integer $n$, we denote by $[n]$ the set $\{1,2,\hdots,n\}$.

\begin{definition}
Let $R$ be a ring. The \emph{$n$th Weyl algebra over $R$} is the $\kk$-algebra $\wanr$ generated over $R$ by $p_1,q_1,\hdots,p_n,q_n$ satisfying the relations:
\begin{align*}
&[q_i,p_i] = 1, [p_i,r]=[q_i,r]=0 & &\text{for all $i \in [n]$, $r \in R$} \\
&[q_i,p_j] = [p_i,p_j] = [q_i,q_j] = 0  & &\text{for all $i,j \in [n]$, $i \neq j$.} 
\end{align*}
\end{definition}

We refer to $\fwa$ as the \emph{first Weyl algebra}\footnote{At some point the brave choice was made to use $\fwa$ instead of the more common $A_1(\kk)$. This is to avoid conflicts later on. We will mostly reserve $A$ for a generic $\kk$-algebra.} and in this setting we write $p=p_1$ and $q=q_1$. This is the canonical example of an infinite-dimensional simple ring. 

What constitutes a proper ``generalization" of the Weyl algebra? Should it be some class of algebras containing the Weyl algebra and whose members have similar algebraic and homological properties? If so, which properties should it maintain and which are we allowed to weaken? Or should it be some class of algebras that are obtained from the Weyl algebra through a formal deformation or twisting procedure?

There is not a correct answer to this question, but in this survey we hope to give an introduction into some of the different branches of study inspired by the Weyl algebra. Of course, it is not possible to go out on every limb. For example, we make almost no mention of $\mathcal{D}$-modules, which some will find such an egregious omission that they immediately don their cap and exit the office\footnote{We make no apologies.}.
The author has made every attempt to be as comprehensive as possible with references, at a cost of a bibliography that is roughly the same size as the paper itself\footnote{Only a mild exaggeration.}.
Undoubtedly, something or someone was forgotten and the author offers sincere apologies for these oversights\footnote{Those responsible have been sacked.}.

\subsection*{Acknowledgements} 
I am indebted to my collaborators, Daniele Rosso and Robert Won, for teaching me key aspects of the theory discussed below. Parts of this article were written while the author was a tag-a-long participant at the 2022 Spannocchia Writer's Conference. The author thanks Tenuta Di Spannocchia for its hospitality and the other conference participants for interesting discussions mostly unrelated to the Weyl algebra. 

\section{Generalized and generalizeder Weyl algebras}

In his 1977 manuscript on Quillen's lemma \cite{joseph}, Joseph introduced a family of subalgebras of $\fwa$ which we recall here\footnote{Joseph denotes these algebras $A_{1,f}$. For reasons, we have taken some liberties with the definition and notation. In particular, we have adjusted the roles of the variables.}. Given $a(u) \in \kk[u]$, the $\kk$-algebras $J(a)$ are generated by $x,y,z$, and satisfy the relations
\[ xz=(z-1)x, \quad yz=(z+1)y, \quad yx=a(z).\]
It is easy to see that $J(a) \iso \fwa$ if and only if $\deg a = 1$. Now if $\deg a > 1$, write $a(u)=(u-\alpha)(u-\beta)b(u)$ for some $b \in \kk[u]$. Making the identification
\[ y=q, \quad z=qp+\alpha, \quad x=-p(-qp-\alpha+\beta) g(-qp-\alpha),\]
shows that $J(a)$ is isomorphic to a subalgebra of $\fwa$. This implies $xyx = xa(z) = a(z-1)x$ and because $\fwa$ is a domain we have that $xy=a(z+1)$. Joseph proves that $J(a)$ is simple if and only if $f$ has no pair of roots $r_1,r_2$ satisfying $r_1-r_2 \in \NN^+$ \cite[Lemma 3.1]{joseph}. We say such a pair of roots are \emph{congruent}. Note that the Weyl algebra itself vacuously satisfies this condition.

These algebras have appeared repeatedly in various guises and generalizations. One consistent connection is with the enveloping algebra of the Lie algebra $\fsl_2$ and its primitive quotients. Recall that $U=U(\fsl_2)$ is the $\kk$-algebra generated by $e,f,h$ with relations
\[ [e,f]=h, \quad [e,h]=2e, \quad [f,h]=-2f.\]
The \emph{Casimir element} of $U$ is $\Omega=4fe+h^2+2h$ and the center of $U$ is $\kk[\Omega]$. The primitive factors are of the form\footnote{The primitive factors $U_{\lambda}$ and $U_{-\lambda}$ are isomorphic. Otherwise these algebras are pairwise distinct.} $U_\lambda=U/(\Omega-(\lambda^2-1))$ for $\lambda \in \kk$. The theory of these algebras owes a huge debt to the work of Dixmier \cite{dix2,dix3}, Smith \cite{sprim}, Stafford \cite{St2}, and others. 

Fix $\lambda \in \kk$ and set $a=\frac{1}{4}(\lambda^2-1)-z^2-z$. Then there is a map $\phi:U \to J(a)$ given by $e \mapsto -y$, $f \mapsto -x$, and $h \mapsto -2z$. It is easy to check that $\phi$ is a surjective homomorphism and moreover, $\ker\phi=(\Omega-(\lambda^2-1))$. It follows\footnote{This is essentially \cite[Example 4.7]{H1}.} that $U_\lambda \iso J(a)$.

In his 1993 paper\footnote{This paper was first received in 1990, which is relevant for the coming discussion.}, Hodges studies the algebras $J(a)$ under the name \emph{deformations of type $A$ Kleinian singularities} \cite{H1}. To make sense of this terminology, we need only note that there is a filtration on $J(a)$ obtained by setting $\deg z = 2$ and $\deg x = \deg y = \deg_u a$. The associated graded ring of this filtration is then isomorphic to $\kk[x,y,z]/(yx-z^{n+1})$, a (commutative) Kleinian singularity. Hodges also notes that these algebras could be viewed as infinite-dimensional primitive factors of a class of algebras studied by Smith \cite{Sm2}. 

But as with many things in mathematics, there seemed to be something ``in the air" regarding these algebras and their generalizations\footnote{The reader is encouraged to remember that this is all before the internet made publicizing ones work trivial.}. In an unpublished 1989 preprint \cite{rosenberg2}, Rosenberg introduced a class of algebras known as \emph{hyperbolic rings}. This work was apparently not widely circulated, but made the rounds at Harvard\footnote{This is pretty much just speculation based on other papers that cite it. If someone has a copy of this paper, please forward to gaddisj@miamioh.edu.}. This construction appeared later in a 1995 manuscript based on lectures spanning from 1991 to 1993 \cite{rosenberg}. Jordan on the other hand was working on certain classes of iterated Ore extensions, called \emph{ambiskew polynomial rings} as well as their primitive factors \cite{Jkrull}. The algebras $J(a)$ can be recovered from this construction. We will return to this below.

\subsection{Generalized Weyl algebras}
The term that has caught on in the literature, however, is \emph{generalized Weyl algebra}, which originally appeared in Bavula's 1990 thesis \cite{B0}. The standard introduction to GWAs is his paper \cite{B1}, published originally in Russian in 1992 and then in English in 1993. 

\begin{definition}\label{def.gwa}
Let $R$ be a $\kk$-algebra, 
$\bsigma=\{\sigma_1,\hdots,\sigma_n\}$ a set of commuting automorphisms of $R$, and $\ba=\{a_1,\hdots,a_n\}$ a set of regular central elements in $R$. The (rank n) \emph{generalized Weyl algebra}\footnote{There is some freedom in the definition. Some require $R$ to be commutative, or even a commutative domain, but we make no such requirements. One can also drop the requirement that the $a_i$ be regular elements of $R$. These are sometimes called \emph{degenerate} GWAs. Those intending to drop the other two requirements are in for a bad time.} (GWA) $R(\bsigma,\ba)$ associated to this data is the $\kk$-algebra generated over $R$ by $x_i,y_i$, $i \in [n]$, satisfying the relations:
\begin{align*}
&x_i r = \sigma_i(r) x_i \quad y_ir = \sigma_i\inv(r) y_i  & &\text{for all $r \in R$ and $i \in [n]$} \\
&y_ix_i = a_i  \quad x_iy_i = \sigma_i(a_i)  & &\text{for all $i \in [n]$} \\
&[x_i,y_j] = [x_i,x_j] = [y_i,y_j] = 0 & &\text{for all $i,j \in [n]$, $i \neq j$}.
\end{align*}
\end{definition}

Let us for the time being consider only the rank one case, so we set $x=x_1$ and $y=y_1$. A rank one GWA $\kk[z](\sigma,a)$ is called \emph{classical} if $\sigma(z)=z-1$ and \emph{quantum}\footnote{Quantum GWAs are also allowed at times to have base ring $\kk[z^{\pm 1}]$, though these behave more like the classical GWAs.} if $\sigma(z)=qz$ for some $q \in \kk^\times$. The classical GWAs are the algebras $J(a)$. Note that any GWA over $\kk[z]$ is, up to isomorphism, either classical or quantum.

A classical GWA $\kk[z](\sigma,a)$ with $a \in \kk^\times$ is isomorphic to $\kk[z][x^{\pm 1};\sigma]$, a skew Laurent ring. When $a$ is linear, $\kk[z](\sigma,a) \iso \fwa$. The reader is encouraged to work out this isomorphism explicitly, at which point it becomes trivial to see how one would obtain the $n$th Weyl algebra from Definition \ref{def.gwa}. By the above discussion, a classical GWA is isomorphic to some $U_\lambda$ when $a$ is quadratic. Many other important rings arise through GWA constructions. This includes $U(\fsl_2)$ itself, the quantum enveloping algebra $U_q(\fsl_2)$ and its primitive factors, as well as the down-up algebras of Benkart and Roby \cite{BRdu}. 

By not requiring $R$ to be commutative, we can define a rank $n$ GWA $R(\bsigma,\ba)$ iteratively. For $k \in [n]$, set $\bsigma_k=\{\sigma_1,\hdots,\sigma_k\}$ and $\ba_k=\{a_1,\hdots,a_k\}$. Then $R(\bsigma_k,\ba_k)$ is the $R$-subalgebra of $R(\bsigma,\ba)$ generated by $x_i,y_i$, $i \in [k]$. Hence, $R(\bsigma_k,\ba_k)$ is a rank $k$ GWA and for $k \in \{2,\hdots,n\}$, we have $R(\bsigma_k,\ba_k)=R(\bsigma_{k-1},\ba_{k-1})(\sigma_k,a_k)$. Thus, $R(\bsigma_n,\ba_n)=R(\bsigma,\ba)$.

Let $R(\bsigma,\ba)$ and $R'(\bsigma',\ba')$ be GWAs of rank $m$ and $n$, respectivesly. Set $\btau=(\tau_1,\hdots,\tau_{m+n})$ and $\bb=(b_1,\hdots,b_{m+n})$ where
\[
\tau_i = \begin{cases}
\sigma_i \tensor 1 & i \leq m \\
1 \tensor \sigma_{i-m}' & m < i \leq n,
\end{cases} \qquad
b_i = \begin{cases}
a_i \tensor 1 & i \leq m \\
1 \tensor a_{i-m}' & m < i \leq n.
\end{cases}
\]
Then $R(\bsigma,\ba) \tensor R'(\bsigma',\ba') \iso (R \tensor R')(\btau,\bb)$ is a rank $m+n$ GWA. Thus, the tensor product of two rank 1 GWAs is a rank 2 GWA, but not every rank 2 GWA is a tensor product of rank 1 GWAs.

\subsection{Quantizations}

Here we describe the quantum Weyl algebras of Giaquinto and Zhang \cite{GZ}.

Let $V$ be an $n$-dimensional vector space and let $\cR:V \tensor V \to V \tensor V$ be a linear transformation. Set $\cR_{12}=\cR \tensor \id_V$ and $\cR_{23} = \id_V \tensor \cR$. Then $\cR$ is a \emph{Hecke symmetry} if
\begin{enumerate}
\item $\cR_{12}\cR_{23}\cR_{12} = \cR_{23}\cR_{12}\cR_{23}$ (the braid relation), and
\item $(\cR-q)(\cR+q\inv) = 0$ for some $q \in \kk^\times$ (the Hecke condition).
\end{enumerate}
Set $\cR(x_i \tensor x_j)=\cR_{ij}^{k\ell}x_k \tensor x_\ell$ with $\cR_{ij}^{k\ell} \in \kk$.

\begin{definition}
Let $\cR$ be a Hecke symmetry. The \emph{quantum Weyl algebra} associated to $\cR$ is the $\kk$-algebra generated by $x_i,y_i$, $i \in [n]$, satisfying the relations:
\[ \sum_{k,\ell} \cR_{ij}^{k\ell} x_k x_\ell - qx_ix_j = 0, \quad
\sum_{k,\ell} \cR_{\ell k}^{ji} y_k y_\ell - qy_iy_j = 0, \quad
y_ix_j = \delta_{ij} + q \sum_{k,\ell} \cR_{j\ell}^{ik} x_k y_\ell, \quad
\text{ for all $i,j \in [n]$}.\]
\end{definition}

One of the most important examples of this construction are the multiparameter quantum Weyl algebras. We recall this definition now, but note that we are using instead the convention of Alev and Dumas \cite{AD2} (see also \cite{Jsimp}), as it is more common in the literature. Another interesting case that arises in \cite{GZ} corresponds to the \emph{Jordan symmetry}.

\begin{definition}
Let $\bq=(q_i) \in (\kk^\times)^n$ and let $\bp=(p_{ij}) \in M_n(\kk^\times)$ be multiplicatively antisymmetric\footnote{Multiplicatively antisymmetric means $p_{ii}=1$ and $p_{ij} = p_{ji}\inv$ for $i \neq j$.}. The (rank $n$) \emph{quantized Weyl algebra} $W_n^{q,\bp}(\kk)$ is generated over $\kk$ by $x_i,y_i$, $i \in [n]$, satisfying the relations:
\begin{align*}
&x_ix_j = q_i p_{ji}\inv x_jx_i, \quad y_jy_i = p_{ji} y_i y_j, \quad
x_iy_j = p_{ji} y_j x_i, \quad x_jy_i = q_i p_{ji}\inv y_ix_j, & &\text{for all $i,j \in [n]$, $i < j$}. \\
&x_iy_i - q_i y_i x_i  = 1 + \sum_{k=1}^{i-1} (q_k-1) y_kx_k, & &\text{for all $i \in [n]$}.
\end{align*}
\end{definition}

We denote the \emph{first quantum Weyl algebra}\footnote{The matrix $\bp$ is irrelevant when $n=1$.} by $W_1^q(\kk)$ (where $q=q_1$) which is generated by $x=x_1$ and $y=y_1$ with relation $xy-qyx=1$. It is not difficult to see that $W_1^q(\kk)$ can be realized as a quantum GWA much the same way that $\fwa$ can be realized as a classical GWA. However, the higher quantized Weyl algebras are not GWAs (at least, not in any obvious way). For this we need a more general construction.

\subsection{TGWAs}\label{sec.TGWA}

It would be reasonable, in light of the above discussion, to hope that enveloping algebras of simple Lie algebras, or at least their primitive quotients, appear as GWAs. In general this does not happen, but a construction of Mazorchuk and Turowksa captures many of these \cite{MT}. Though they did not make this explicit connection, they did observe an analogue between certain module categories.

\begin{definition}
Let $R$ be a $\kk$-algebra, $\sigma=\{\sigma_1,\hdots,\sigma_n\}$ a set of commuting automorphisms of $R$, $a=\{a_1,\hdots,a_n\} \subset R$ a set of central nonzero divisors in $R$, and $\mu=(\mu_{ij}) \in M_n(\kk^\times)$. The \emph{twisted generalized Weyl construction (TGWC)} is the $\kk$-algebra $\cC_\mu(R,\sigma,a)$ generated over $R$ by $X_i^{\pm}$, $i \in [n]$ satisfying the relations: 
\[
X_i^{\pm} r - \sigma_i^{\pm 1}(r) X_i^{\pm}, \quad
X_i^-X_i^+ - a_i, \quad 
X_i^+X_i^- - \sigma_i(a_i), \quad
X_i^+ X_j^- - \mu_{ij} X_j^- X_i^+, \quad \text{for all $i,j \in [n], i \neq j, r \in R$}.
\]
There is a $\ZZ^n$-grading\footnote{See Section \ref{sec.grading}.} on $\cC_\mu(R,\sigma,a)$ is obtained by setting $\deg(R)=\bzero$ and $\deg(X_i^{\pm})=\pm \be_i$ for all $i$. 
The associated \emph{twisted generalized Weyl algebra (TGWA)} is the quotient $\cA_\mu(R,\sigma,a) = \cC_\mu(R,\sigma,a)/\cJ$ where $\cJ$ is the sum of all $\ZZ^n$-graded ideals $J$ with $J_\bzero = \{0\}$. 
\end{definition}

We recover the class of rank $n$ GWAs from those TGWAs with $\mu_{ij}=1$ for all $i,j$ and $\sigma_i(a_j) = a_j$ for $i \neq j$. Quantized Weyl algebras \cite{FH1}, certain primitive quotients of $U(\fsl_{n+1})$ and $U(\fsp_{2n})$ \cite{HVS}, and certain primitive quotients of enveloping algebras of Lie superalgebras \cite{HVS2} can all be realized as TGWAs. Under mild restrictions, TGWAs are closed under tensor products in a similar way as GWAs \cite{GR}. An entire survey could be written just about TGWAs\footnote{I smell a sequel!}, so we can only touch on some of the recent research in this area.

We will pay special attention to a certain subclass of TGWAs. Suppose $\cA=\cA_\mu(R,\sigma,a)$ is a TGWA of rank $n$ such that each $a_i$ is a nonzero divisor in $R$ ($\cA$ is \emph{regular}) and that the canonical map $R \to \cA$ is injective ($\cA$ is \emph{$\mu$-consistent}). For each $i,j \in [n]$ with $i<j$, set 
$V_{ij} = \Span_\kk\{\sigma_i^k(t_j) : k \in \ZZ\}$.
If each $V_{ij}$ is finite-dimensional over $\kk$, then we say that $\cA$ is \emph{$\kk$-finitistic} \cite{hart3}. 

We specialize just a bit more. Suppose $\sigma_i(t_j) = \gamma_{ij} t_j$ for all $i,j \in [n]$ with $i \neq j$. In this case, $\cA$ is consistent if and only if $\mu_{ij}\mu_{ji} = \gamma_{ij}\gamma_{ji}$ for all $i \neq j$. Then $\cA$ has the usual TGWA relations, and the ideal $J$ is generated by the elements:
\[ 
X_i^+ X_j^+ - \gamma_{ij}\mu_{ij}\inv X_j^+ X_i^+, \quad
X_j^- X_i^- - \gamma_{ij} \mu_{ji}\inv X_i^- X_j^-, \qquad i,j \in [n], i \neq j.
\]
These are known as TGWAs of type $(A_1)^n$ and the GWAs of rank $n$ belong to this subclass. The work of Futorny and Hartwig \cite{FH1} gives an explicit connection between these TGWAs and quantized Weyl algebras.

\subsection{BR algebras}\label{sec.BR}

In their classification of simple $\ZZ$-graded rings of Gelfand-Kirillov dimension\footnote{See Section \ref{sec.GK}.} two, Bell\footnote{To be clear, this is Jason Bell, which would also be my name would be if my advisor adopted me.} and Rogalski introduced a construction that generalizes rank one GWAs \cite{BR}. One could reasonably think about this definition as a noncommutative extended Rees ring construction. Indeed, these algebras have their origin in the study of noncommutative blowing up \cite{KRS}. The definition below is taken from \cite{GRW1}, which is slightly less restrictive that the original. 

\begin{definition}\label{defn.BR}
Let $R$ be a $\kk$-algebra and let $\sigma$ be an automorphism of $R$. Let $H$ and $J$ be two-sided ideals of $R$. Set $I^{(0)} = R$, $I^{(n)} = J \sigma(J) \cdots \sigma^{n-1}(J)$ for $n \geq 1$ and $I^{(n)} = \sigma^{-1}(H) \sigma^{-2}(H) \cdots \sigma^n(H)$ for $n \leq -1$. We assume that the ideals $H$ and $J$ are chosen so that $I^{(n)} \neq 0$ for all $n$.
The corresponding \emph{Bell--Rogalski algebra}\footnote{Bell and Rogalski didn't name these algebras. At one point I suggested that they should be called \emph{generalized generalized Weyl algebras} but my coauthors rejected it. I am now in favor of \emph{generalizeder Weyl algebras}.} (or \emph{BR algebra}, for short) is the algebra
\[ R(t,\sigma,H,J) = \bigoplus_{n \in \ZZ} I^{(n)} t^n \subseteq R[t,t^{-1};\sigma].\]
\end{definition}

One can obtain the GWA $R(\sigma,a)$ above by taking $H=R$ and $J=(a)$. The identification comes from the map $y \mapsto t\inv$ and $x \mapsto at$. Conversely, every BR algebra with $H$ and $J$ principal ideals is isomorphic to a GWA \cite{GRW1}. We will see that many of the properties we can study for GWAs can actually be considered in this more general context.

\section{Ring-theoretic properties}

Below we discuss some of the ways that GWAs generalize properties of the Weyl algebra, but more precisely those of primitive quotients of $U(\fsl_2)$. When possible, we also discuss these properties in the context of other generalizations introduced above.

\subsection{Grading}\label{sec.grading}

Let $\Gamma$ be an abelian group. A $\kk$-algebra $A$ is $\Gamma$-graded if there exists a vector space decomposition $A=\bigoplus_{\gamma \in \Gamma} A_{\gamma}$ such that $A_{\gamma} A_{\gamma'} \subset A_{\gamma+\gamma'}$. A left $A$-module $M$ is $\Gamma$-graded if it has its own vector decomposition $M=\bigoplus_{\gamma \in \Gamma} M_{\gamma}$ such that $A_{\gamma} M_{\gamma'} \subset M_{\gamma+\gamma'}$. We denote by $\GrMod_{\Gamma} A$ the category of graded left $A$-modules along with degree zero homomorphisms.

The Weyl algebra $\wanr$ over a ring $R$ has a canonical $\ZZ^n$-grading obtained by setting $\deg(p_i)=\be_i$ and $\deg(q_i)=-\be_i$ for all $i \in [n]$, and $\deg(r)=0$ for all $r \in R$. Similarly, a TGWA $\cA_\mu(R,\sigma,a)$ inherits its $\ZZ^n$-grading from the corresponding TGWC, so $\deg(X_i^{\pm})=\pm \be_i$ for all $i \in [n]$ and $\deg(r)=\bzero$ for all $r \in R$. There is a natural $\ZZ$-grading on a BR algebra inherited from $R[t^{\pm 1};\sigma]$, where $\deg(r)=0$ and $\deg(t^{\pm 1}) = \pm 1$.

\subsection{Domains}
We order the monomials of the $R$-basis of $R(\sigma,a)$ degree-lexicographically (with $y>x$). Let $f_1,f_2 \in R(\sigma,a)$ be nonzero and choose the monomial with maximal degree in each, say $m_1=x^{i_1} y^{j_1}$ and $m_2 = x^{i_2} y^{j_2}$, respectively. An inductive argument shows $m_1m_2 = x^{i_1+i_2} y^{i_1+i_2} + (\text{lower degree terms})$. Hence, the monomial of maximal degree in $A$ is $x^{i_1+i_2} y^{i_1+i_2}$, so $f_1f_2 \neq 0$. It follows that if $R$ is a domain, then so is $R(\sigma,a)$. By using the iterative construction, this implies that the property of being a domain extends to higher rank GWAs $R(\bsigma,\ba)$ assuming $R$ is a domain.

This is one property that holds for all generalizations of the Weyl algebra mentioned above, though for quantizations additional assumptions are needed (see \cite{GZ}). For TGWAs see \cite[Proposition 2.9]{FH1}, and for BR algebras see \cite[Lemma 3.1]{GRW1}.

\subsection{Simplicity}
A ring is said to be \emph{simple} if its only two-sided ideals are the zero ideal and the ring itself. Since we have assumed $\kk$ has characteristic zero, $\wan$ is simple. A standard proof of this for $n=1$ takes a nonzero element $b \in I$ and uses repeated applications of commuting with $p$ and $q$ to show that $I$ contains a nonzero constant\footnote{Alternatively, observe that a ring is simple if it has no nonzero finite-dimensional representation. This is equivalent to the given definition because the annihilator of a finite-dimensional representation is necessarily a proper nonzero two-sided ideal. Now let $X,Y \in M_n(\kk)$ be an $n$-dimensional representation of $\fwa$, so $XY-YX=I_n$. But then $0=\trace(XY-YX)= \trace(I_n)=n$, a contradiction.}. The proof for higher $n$ is similar. 
Over a field $F$ of positive characteristic, the Weyl algebra $W_n(F)$ is not simple.

There is interest in understanding infinite-dimensional simple rings in a broader context. By the result of Joseph mentioned earlier, a classical GWA $\kk[z](\sigma,a)$ is simple precisely when $a$ has no pair of congruent roots\footnote{Repeated roots are ok.}. There is, however, a more general criteria for rank one GWAs. This can be attributed both to to Bavula \cite[Theorem 4.2]{B3} and Jordan \cite[Theorem 6.1]{Jprim}.

An ideal of $R$ is \emph{$\sigma$-stable} if $\sigma(I)=I$. The automorphism $\sigma$ is said to be \emph{inner} if there exists a unit $\alpha \in R$ such that $\sigma^n(r)=\alpha\inv r \alpha$ for all $r \in R$. Then $R(\sigma,a)$ is simple if and only if $R$ has no proper $\sigma$-stable ideal, no power of $\sigma$ is an inner automorphism of $R$, and $R=Ra+R\sigma^n(a)$ for all $n \geq 1$. This last condition is equivalent to the congruent roots condition for classical GWAs.

A general simplicity criterion for TGWAs was given by Hartwig and \"{O}inert \cite{HO}. Here we will only recall their result for those of type $(A_1)^n$. Let $A=A_\mu(R,\sigma,t)$ be such a TGWA. Then $A$ is simple if and only if $R\sigma_i^k(t_i)+Rt_i=R$ for all $i \in [n]$ and all $k \in \ZZ_{>0}$, $R$ is $\ZZ^n$-simple\footnote{Here $\ZZ^n$-simplicity is nothing more than simplicity with respect to the group generated by $\sigma$.}, and $\cnt(A) \subset R$. 

For a BR algebra $B=R(t,\sigma,H,J)$, set 
\begin{align}\label{eq.SB}
\cS(B) = \{\bp\in\Spec(R) : \bp\supset HJ\}=\Spec(R/(HJ))=\Spec(R/H)\cup \Spec(R/J).
\end{align}
We say $\cS(B)$ is \emph{$\sigma$-lonely} if $\sigma^i(\cS(B)) \cap \cS(B) = \emptyset$ for all $i \neq 0$. Then $B$ is simple if and only if $R$ is $\sigma$-simple and $\cS(B)$ is $\sigma$-lonely \cite[Proposition 2.18 (3)]{BR} and \cite[Proposition 2.3]{GRW1}. It is not difficult to see that this criterion reduces the GWA criterion above.

For quantizations, there is no general criteria for simplicity, though certainly they can be, e.g. the higher Weyl algebras themselves, see also \cite[Theorem 3.5]{GZ}. The higher quantized algebras are not simple in general. In this case, however, Jordan has presented a localization which is simple \cite{Jsimp}, a result which will appear again later in this survey.

\subsection{Noetherian}
We begin with a construction that is useful in studying properties of GWAs. A \emph{$\sigma$-derivation} of $R$ is a $\kk$-linear map $\delta$ satisfying $\delta(rr') = \sigma(r)\delta(r') + \delta(r)r'$ for all $r,r' \in R$. Given a $\sigma$-derivation $\delta$ of $R$, the \emph{Ore extension} $R[x;\sigma,\delta]$ is generated over $R$ by $x$ satisfying the rule $xr=\sigma(r)x + \delta(r)$ for all $r \in R$.

Consider the Ore extension $T=R[u;\sigma][v;\sigma\inv,\delta]$ where $\sigma$ is extended to $R[u;\sigma]$ by setting $\sigma(u)=u$ and $\delta$ is determined by $\delta(r)=0$, $\delta(u)=a-\sigma(a)$ for some $a \in R$. This iterated extension is what Jordan calls an \emph{ambiskew polynomial ring} \cite{Jiter}.
One observes that $vu-a$ is a central element in $T$. Thus, $vu-a$ generates a two-sided ideal of $T$ and it is clear that $T/(vu-a) \iso R(\sigma,a)$.

Now suppose $R$ is (left or right) noetherian. Then $T$ is noetherian by the Hilbert Basis Theorem for Ore extensions, see \cite{GW}. Since $T/(vu-a) \iso R(\sigma,a)$, then the noetherian property passes down to $R(\sigma,a)$. Another proof of this fact uses the Hilbert Basis Theorem applied directly to GWAs, see \cite[Proposition 1.3(1)]{B1}. As with the property of being a domain, higher rank GWAs over a noetherian ring are noetherian.

The noetherian property for other generalizations is more subtle. By \cite[Theorem 4.11(1)]{GR}, TGWAs of type $(A_1)^n$ over a noetherian ring are noetherian. It is also true that certain other TGWAs over noetherian rings are noetherian, but there is no general argument of this nature.

For a BR algebra $B=R(t,\sigma,H,J)$, recall the set $\cS(B)$ defined in \eqref{eq.SB}. We say that $\cS(B)$ is \emph{critically dense} in $\Spec R$ if for any proper closed subset $Y \subset X$, $Y \cap \bp = \emptyset$ for all but finitely many $\bp \in \cS(B)$. If $\cS(B)$ is critically dense, then $B$ is noetherian \cite[Proposition 2.18 (1)]{BR}. However, critical density is not equivalent to the noetherian condition, see \cite{GRW2}.

\subsection{GK dimension}\label{sec.GK}
Here we only recall the definition sufficiently for our setting. Let $A$ be an affine $\kk$-algebra with generating subspace $V$. Set $A_m = \sum_{k=0}^m V_k$ and $d_m = \dim_\kk A_m$. The \emph{Gelfand-Kirillov (GK) dimension} of $A$ is
\[ \GKdim A = \limsup \frac{ \log( \dim_k A_m)}{\log m}.\]
For a thorough overview of GK dimension we refer the reader to \cite{KL}.
It is worth noting that GK dimension is preserved upon taking associated graded rings, and that for commutative rings GK dimension is equal to Krull dimension. Under the standard filtration, $\gr \fwa \iso \kk[p,q]$, so $\GKdim \fwa = 2$. Similarly, using the filtration discussed above, any classical GWA has GK dimension 2.

There is a more general result, however, which was proved almost simultaneously by Ebrahim \cite[Theorem 2.6]{ebrahim} and Zhao-Mo-Zhang \cite[Theorem 3.4]{ZMZ}. We say the automorphism $\sigma$ of $R$ is \emph{locally algebraic}\footnote{Equivalently, $\dim_\kk \{\sigma^n(r) : n \geq 0\} < \infty$ for all $r \in R$.} if any finite-dimensional subspace $U$ of $R$ is contained in a $\sigma$-stable finite subspace of $R$. Note that if $\sigma$ is finite order, then this condition is automatic. If $\sigma$ is a locally algebraic automorphism of a $\kk$-algebra $R$, then $\GKdim R(\sigma,a) = \GKdim R + 1$. One can also apply the GK dimension argument iteratively to higher rank GWAs assuming that each $\sigma_i$ is locally algebraic.

This condition is useful for determining the GK dimension of a BR algebra $B=R(t,\sigma,H,J)$ as well. Suppose $\kk$ is algebraically closed, let $K$ be the field of fractions of $R$, and set $d=\trdeg(K/\kk)$. If $\sigma$ is locally algebraic, then $\GKdim B = d +1$.

It is likely that one could obtain a similar result in the case of TGWAs of type $(A_1)^n$, but no such result exists in general in the literature. A quantized Weyl algebras of rank $n$, which the reader will recall is an example of such an algebra, has GK dimension $2n$ \cite[Proposition 3.4]{GL3}. 

\subsection{Global dimension}
The \emph{left (resp. right) global dimension} of a ring $A$ is the supremum of the projective dimensions of left (resp. right) $A$-modules $M$. For a noetherian ring, left and right global dimension coincide. Thus, global dimension measures how close a ring is to being semisimple artinian.

When $\chr\kk=0$, $\gldim \fwa = 1$, that is, $\fwa$ is hereditary, see \cite[Chapter 7]{MR}. For the $U_\lambda$, we have a trichotomy: $\gldim U_0=\infty$, $\gldim U_m = 2$ for nonzero integers $m$, and $\gldim U_\lambda = 1$ otherwise. This is due to Stafford \cite{St2} and was motivated by earlier results of Roos \cite{roos1} and Smith \cite{sprim}. Viewing $U_\lambda$ as a GWA, we have $a=\frac{1}{4}(\lambda^2-1)-z^2-z$. Hence, the condition $\lambda=0$ implies that $a= -(z+1/2)^2$, so $a$ has a repeated root. One the other hand, the condition $\lambda \in \ZZ\backslash\{0\}$ means that the two roots of $a$ differ by an integer.

This leads to a general criteria for rank one GWAs. A classical GWA $\kk[z](\sigma,a)$ has infinite global dimension when $a$ has a repeated root, global dimension two when $a$ has no repeated roots but a pair of congruent roots, and global dimension one otherwise. This result is due to Hodges \cite[Theorem 4.4]{H1}. For a more general argument, see \cite[Corollary 7.8]{Jkrull}. This extends to GWAs $R(\sigma,a)$ over a Dedekind ring $R$ where the role of roots is replaced by ideals $I$ in $R$ such that $a \in I$. See also \cite{bavgldim}.

A quantized Weyl algebra of rank $n$ has global dimension $2n$ except in special cases (see \cite[Theorem 3.5]{GZ} and \cite[Theorem 3.8]{FKK}). There is no general criteria for global dimension of BR algebras or TGWAs. See the work of McConnell for the case of Weyl algebras over affine rings \cite{Mc}, and the work of Goodearl, Hodges, and Lenagan for results in the case of Weyl algebras over division rings \cite{GHL}.

The Krull dimension of a GWA has been studied alongside that of global dimension in many studies. We do not review all of these results, but point the interested reader to handful of handy references \cite{BL2,BL1,BVO,FKK,GHL,Jkrull,Mc}.

\section{Representation theory and other modular considerations}

The classification of irreducible representations of the first Weyl algebra $\fwa$, as well as those of $U(\fsl_2)$, is due to Block \cite{block}. The representation theory of GWAs draws heavily from that of Lie algebras (and their enveloping algebras). Throughout this section, we assume $R$ is commutative and let $\Maxspec(R)$ denote the set of maximal ideals of $R$.

\subsection{Weight modules}

We begin by describing some (left) modules of the first Weyl algebra $\fwa$. Let $V^+=\{v_i\}_{i \in \ZZ_{\geq0}}$ and define an $\fwa$-action by $pv_i=v_{i+1}$ and $qv_i = i v_{i-1}$ for all $i$. Similarly, let $V^- = \{v_i\}_{i \in \ZZ_{<0}}$ and define an $\fwa$-action by $qv_i = i v_{i-1}$ for all $i$, $pv_{-1}=0$, and $pv_i=v_{i+1}$ for $i < -1$. Finally, for $\lambda \in (0,1)$, let $V^\lambda = \{ v_i \}_{i \in \ZZ}$ and define an $\fwa$-action by $pv_i=v_{i+1}$ and $qv_i = (i+\lambda) v_{i-1}$ for all $i$. It is easy to check that these actions are well-defined and that each is an (infinite-dimensional) simple module of $\fwa$. Moreover, the modules are non-isomorphic. We now generalize this discussion to modules for rank one GWAs.

We say a (left) module $M$ of the GWA $R(\sigma,a)$ is a \emph{weight module} if there is a decomposition 
\[ M = \bigoplus_{\fm \in \Maxspec(R)} M_{\fm}\]
where $M_{\fm} = \{ z \in M : \fm \cdot z = 0\}$. We call $M_{\fm}$ a \emph{weight space} of $M$ and the \emph{support} of $M$ is
\[ \supp(M) = \{ \fm \in \Maxspec(R) : M_{\fm} \neq 0\}.\]
The analogy here is that $R$ plays the role of the Cartan subalgebra of a Lie algebra.

Suppose $M$ is a weight module for $R(\sigma,a)$. There is an action of $\ZZ$ on $\Maxspec(R)$ given by
$n \cdot \fm = \sigma^n(\fm)$. Then $x(M_{\fm}) \subset M_{\sigma(\fm)}$ and $y(M_{\fm}) \subset M_{\sigma\inv(\fm)}$. In fact, $M$ splits as 
\[ M = \bigoplus_{\cO \in \Maxspec(R)/\ZZ} M_{\cO} \quad\text{where}\quad \supp(M_{\cO}) \subset \cO.\]
Thus, in classifying weight modules for a particular GWA, it suffices to consider those that are supported on one orbit at a time. 
Bavula gave a classification of simple and indecomposable weight modules for classical GWAs \cite{B1}. However, the approach that we follow here is the more closely related to the work of Drozd, Guzner, and Ovzienko \cite{DGO}.

We say an ideal $\fm \in \Maxspec(R)$ is a \emph{break} if $a \in \fm$. In the case of $R=\kk[t]$, the breaks correspond to the roots of the defining polynomial $a$. If $\fm$ is break, then it is straightforward to show that $xM_{\fm}=yM_{\sigma(\fm)}=0$. Thus, the simple weight modules correspond to intervals between breaks. The weight spaces are 1-dimensional and isomorphic to $R/\sigma^i(\fm)$. We can visualize as follows:

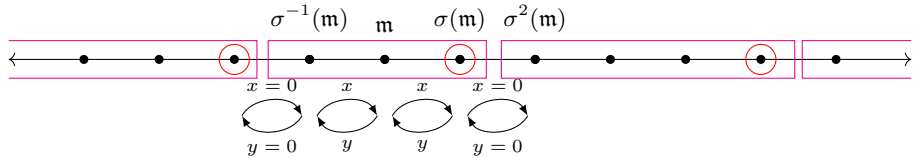
\begin{center}
\begin{tikzpicture}
\draw[<->] (-6,0) -- (6,0); % x-axis

\fill[black] (-5,0) circle (0.6 mm);
\fill[black] (-4,0) circle (0.6 mm);
\fill[black] (-3,0) circle (0.6 mm);
\fill[black] (-2,0) circle (0.6 mm);
\fill[black] (-1,0) circle (0.6 mm);
\fill[black] (0,0) circle (0.6 mm);
\fill[black] (1,0) circle (0.6 mm);
\fill[black] (2,0) circle (0.6 mm); 
\fill[black] (3,0) circle (0.6 mm);
\fill[black] (4,0) circle (0.6 mm);
\fill[black] (5,0) circle (0.6 mm);

\fill[black] (-2,0) circle (0.6 mm) node[above=7pt] {\small $\sigma\inv(\fm)$};
 \fill[black] (-1,0) circle (0.6 mm) node[above=7pt] {\small $\fm$};
\fill[black] (0,0) circle (0.6 mm) node[above=7pt] {\small $\sigma(\fm)$};
\fill[black] (1,0) circle (0.6 mm) node[above=7pt] {\small $\sigma^2(\fm)$};

\draw[color=red] (-3,0) circle (2mm);
\draw[color=red] (0,0) circle (2mm); 
\draw[color=red] (4,0) circle (2mm); 

\draw[-latex] (.1,-.75) to[out=60,in=120] node[midway,font=\scriptsize,above] {$x=0$} (.9,-.75);
\draw[-latex] (.9,-.75) to[out=-120,in=-60] node[midway,font=\scriptsize,below] {$y=0$} (.1,-.75);
\draw[-latex] (-.9,-.75) to[out=60,in=120] node[midway,font=\scriptsize,above] {$x$} (-.1,-.75);
\draw[-latex] (-.1,-.75) to[out=-120,in=-60] node[midway,font=\scriptsize,below] {$y$} (-.9,-.75);
\draw[-latex] (-1.9,-.75) to[out=60,in=120] node[midway,font=\scriptsize,above] {$x$} (-1.1,-.75);
\draw[-latex] (-1.1,-.75) to[out=-120,in=-60] node[midway,font=\scriptsize,below] {$y$} (-1.9,-.75);
\draw[-latex] (-2.9,-.75) to[out=60,in=120] node[midway,font=\scriptsize,above] {$x=0$} (-2.1,-.75);
\draw[-latex] (-2.1,-.75) to[out=-120,in=-60] node[midway,font=\scriptsize,below] {$y=0$} (-2.9,-.75);

\draw[color=magenta] (-6,.25)--(-6+3.3,.25)--(-6+3.3,-.25)--(-6,-.25);
\draw[color=magenta] (-2.55,.25) rectangle ++ (2.9,-.5);
\draw[color=magenta] (0.55,.25) rectangle ++ (3.9,-.5);
\draw[color=magenta] (6,-.25)--(4.55,-.25)--(4.55,.25)--(6,.25);
\end{tikzpicture}
\end{center}
For the Weyl algebra, there is one break. On the orbit with the break, there are exactly two (infinite-dimensional) simple weight modules. On the other hand, for an orbit without a break, there is exactly one (again, infinite-dimensional) simple weight module.

The indecomposable weight modules can be formed by ``patching" together simple weight modules along breaks. Each break is labeled with and $x$ or a $y$ to indicate which choice we are making as to which one acts nontrivially.

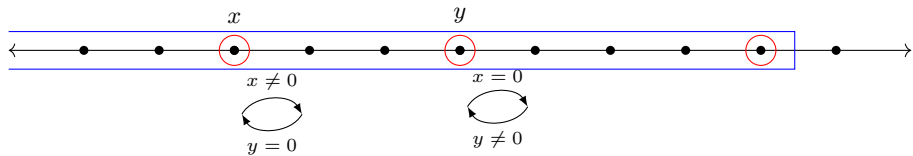
\begin{center}
\begin{tikzpicture}

\draw[<->] (-6,0) -- (6,0); % x-axis
\fill[black] (-5,0) circle (0.6 mm);
\fill[black] (-4,0) circle (0.6 mm);
\fill[black] (-3,0) circle (0.6 mm);
\fill[black] (-2,0) circle (0.6 mm);
\fill[black] (-1,0) circle (0.6 mm);
\fill[black] (0,0) circle (0.6 mm);
\fill[black] (1,0) circle (0.6 mm);
\fill[black] (2,0) circle (0.6 mm); 
\fill[black] (3,0) circle (0.6 mm);
\fill[black] (4,0) circle (0.6 mm);
\fill[black] (5,0) circle (0.6 mm);
\draw[color=red] (-3,0) circle (2mm);
\draw[color=red] (0,0) circle (2mm); 
\draw[color=red] (4,0) circle (2mm);

\fill[black] (-3,0) circle (0.6 mm) node[above=7pt] {\small $x$};
\draw[-latex] (-2.9,-.85) to[out=60,in=120] node[midway,font=\scriptsize,above] {$x\neq 0$} (-2.1,-.85);
\draw[-latex] (-2.1,-.85) to[out=-120,in=-60] node[midway,font=\scriptsize,below] {$y=0$} (-2.9,-.85);
\fill[black] (0,0) circle (0.6 mm) node[above=7pt] {\small $y$};
\draw[-latex] (.1,-.75) to[out=60,in=120] node[midway,font=\scriptsize,above] {$x=0$} (.9,-.75);
\draw[-latex] (.9,-.75) to[out=-120,in=-60] node[midway,font=\scriptsize,below] {$y\neq0$} (.1,-.75);

\draw[color=blue] (-6,.25)--(3.9+.55,.25)--(3.9+.55,-.25)--(-6,-.25);
\end{tikzpicture}
\end{center}

One can also describe simple and indecomposable weight modules in the case of a finite orbit, but we do not do that here. The interested reader is directed to \cite{DGO}.

The description of weight modules for BR algebras is entirely analogous to the above in the case of an infinite orbit \cite{GRW1}. The key difference is that a \emph{break} in defined to be an ideal $\fm \in \Maxspec(R)$ such that $\sigma(m) \in \cS(B)$ (see \eqref{eq.SB}).

It is natural to see how the study of weight modules would extend to higher rank GWAs \cite{B1,BB}. Here the orbits are determined by the action of the group generated by the $\sigma_i$. One can also do this in the case of TGWAs and this has been a driving force behind their development, dating back to the original paper of Mazorchuk and Turowska \cite{MT}. We do not describe this here, but refer the reader to a handful of papers developing the theory \cite{hart1,hart4,HR2}.

\subsection{Morita equivalence}

In \cite{Sm3}, Smith gave an example of a right ideal $I$ of $\fwa$ such that $S=\End_{\fwa}(I)$ is not isomorphic to $\fwa$, but $M_2(S) \iso M_2(\fwa)$. This inspired Stafford to study endomorphisms of right ideals of $\fwa$ more generally. He showed that the above behavior is satisfied for any noncyclic, projective right ideal of $\fwa$ \cite{StEndo}. Stafford went on to classify domains Morita equivalent to $\fwa$, which turn out to be in natural bijection with certain orbits of $\Aut_{\kk}(\fwa)$. Following further work by Musson \cite{musson2}, Berest-Wilson \cite{BW}, and Kouakou \cite{kouakou}, it was established that the isomorphism classes of domains Morita equivalent to $\fwa$ are in bijection with the nonnegative integers. See the work of Sierra for rings \emph{graded} equivalent to $\fwa$ \cite{sierra}.
Such a classification is not known in the setting of GWAs.

Let $\lambda,\mu \in \kk\backslash\{0\}$. By a result of Stafford, $U_\lambda$ is Morita equivalent to $U_{\mu}$ when $\lambda-\mu \in \ZZ$ \cite[Corollary 3.3]{St2}. The converse was proved by Hodges \cite[Theorem 5]{H2}. Recall that $U_0$ is the singular case with infinite global dimension. Equivalently, the defining polynomial in the corresponding GWA has a repeated root. 

Morita equivalence for rank one GWAs was considered by Hodges \cite{H1} and Jordan \cite{Jkrull}. We describe the result only in the classical case. Suppose $a(z) \in \kk[z]$ factors as $a(z)=b(z)c(z)$. Let $a'(z)=b(z-1)c(z)$. The classical GWAs $\kk[z](\sigma,a)$ is Morita equivalent to $\kk[z](\sigma,a')$ \emph{unless} shifting created (or destroyed) a multiple root. The general rank one GWA case is similar, but one considers more generally orbits by the automorphism $\sigma$. See the work of Richard and Solotar for an extension to higher rank (classical) GWAs \cite{RS3}.

\subsection{Graded modules}

For this section we will deal only with $\ZZ$-graded algebras, so set $\GrMod A = \GrMod_{\ZZ} A$.

In some sense, the picture for the graded simple modules is identical to that of the simple weight modules above. For a rank one GWA $R(\sigma,a)$, the key difference is that our choice above of where to place $\fm$ on the number line was arbitrary. In the category $\GrMod R(\sigma,a)$, the choices of where to place $\fm$ reflect graded shifts of the module. Consequently, the classification of simple weight modules also classifies simple graded modules \emph{up to shift}.

For example, in the case of $\fwa$, there are two graded simple modules for each each integer, and one graded simple for each $\lambda \in \kk\backslash\ZZ$ (see Sierra \cite{sierra}). Thus the simples can be realized as the affine line $\kk$ but with each integer point replaced by a double point\footnote{The code for this diagram was stolen shamelessly from a paper of Robert Won without permission.}:

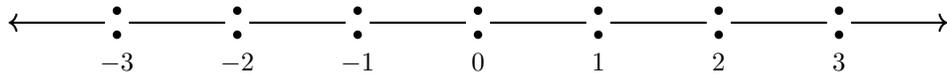
\begin{center}
\begin{tikzpicture}[ scale=1.6]	
			\node at (-1.7,.1)[label = above: {}]{};
		\draw[thick,<-](-3.9,0)--(-3.1,0);
		\fill[black](-3,-.1) circle (1pt);
		\fill[black](-3,.1) circle (1pt);
		\node at (-3,-.1)[label = below: $-3$]{};
		\draw[thick](-2.9,0)--(-2.1,0);
		\fill[black](-2,-.1) circle (1pt);
		\fill[black](-2,.1) circle (1pt);
		\node at (-2,-.1)[label = below: $-2$]{};
		\draw[thick](-1.9,0)--(-1.1,0);
		\fill[black](-1,-.1) circle (1pt);
		\fill[black](-1,.1) circle (1pt);
		\node at (-1,-.1)[label = below:$-1$]{};
		\draw[thick](-.9,0)--(-.1,0);
		\fill[black](0,-.1) circle (1pt);
		\fill[black](0,.1) circle (1pt);
		\node at (0,-.1)[label = below: $0$]{};
		\draw[thick](.1,0)--(.9,0);
		\fill[black](1,-.1) circle (1pt);
		\fill[black](1,.1) circle (1pt);
		\node at (1,-.1)[label = below: $1$]{};
		\draw[thick](1.1,0)--(1.9,0);		
		\fill[black](2,-.1) circle (1pt);
		\fill[black](2,.1) circle (1pt);
		\node at (2,-.1)[label = below: $2$]{};
		\draw[thick,](2.1,0)--(2.9,0);		
		\fill[black](3,-.1) circle (1pt);
		\fill[black](3,.1) circle (1pt);
		\node at (3,-.1)[label = below: $3$]{};
		\draw[thick,->](3.1,0)--(3.9,0);		
\end{tikzpicture}
\end{center}

The general case of a classical GWA $\kk[z](\sigma,a)$ follows similarly. The reader is directed to \cite{bavsimp} and also \cite{woncomm}.

Let $\Tors A$ denote the full subcategory of $\GrMod A$ consisting of sums of finite-dimensional modules and set $\QGrMod A := \GrMod A/\Tors A$. One can also restrict to noetherian objects, for which we use the notation $\qgrmod A := \grmod A/\tors A$. Obviously for a simple GWA $A$, $\QGrMod A = \GrMod A$. As above, suppose $a(z) \in \kk[z]$ factors as $a(z)=b(z)c(z)$. Let $a'(z)=b(z-1)c(z)$. If $A=\kk[z](\sigma,a)$ and $A'=\kk[z](\sigma,a')$, then $\QGrMod A \equiv \QGrMod A'$. This is true \emph{regardless} of whether a multiple root has been created or destroyed. Consequently, every classical GWA is $\QGrMod$-equivalent to a simple GWA. This was proved by Won for the $U_\lambda$ \cite[Theorem 4.19]{woncomm}. See \cite[Theorem 3.8]{FGW} for the general case.

The importance of the categories $\QGrMod$ and $\qgrmod$ come from the study of noncommutative algebraic geometry \cite{AZ}. Temporarily, let $A$ be a connected $\NN$-graded domain of GK dimension 2 that is generated in degree 1. A theorem of Artin and Stafford \cite{ArSt} proves that there exists a projective curve $X$ such that $\qgrmod A \equiv \coh(X)$, the category of coherent sheaves on $X$.

Smith showed that there exists a quotient stack $\cX$ such that $\GrMod\fwa \equiv \Qcoh \cX$ \cite{smith}. We do not describe the quotient stack here explicitly, but note that it comes from a commutative ring built from the Picard group\footnote{The Picard group is the set of autoequivalences modulo natural isomorphisms.} of $\GrMod \fwa$. In \cite{wonpic}, Won computed the Picard group of the $U_\lambda$ use this to give an analog of Smith's result in this case \cite{woncomm}. See \cite{FGW} for Smith's result in the case of a classical GWA.

\section{Invariant Theory}

Classically, invariant theory refers to actions of groups on polynomial rings. Modern work in noncommutative invariant theory has focused on group and Hopf algebra actions on Artin-Schelter regular rings, which may be regarded as noncommutative versions of polynomial rings. We refer the interested reader to the survey of Kirkman for an overview of work in this area \cite{Ki}. In this short section we review some invariant-theoretic results regarding GWAs. Throughout this section we assume that $\kk$ is algebraically closed and $\chr\kk=0$.

\subsection{Automorphisms and Isomorphisms}
A classical problem for any algebra (or algebraic object) is to compute its automorphism group\footnote{More generally, to determine its symmetries.}. Related to this is the isomorphism problem, which asks under what conditions two objects within a class of algebras are isomorphic.

Dixmier computed generators for the automorphism group of $\fwa$ \cite{dix}, as well as those of the algebras $U_\lambda$ \cite{dix3}. This was extended to find generators for the automorphism group of a classical GWA by Bavula and Jordan \cite{BJ}, which we describe here. 

Let $\kk[z](\sigma,a)$ be a classical GWA with $n = \deg_z(a)$. Let $\lambda \in \kk$, $\beta \in \kk^\times$, $m \in \NN$, and let $\Delta_m$ be the linear map $\kk[z] \rightarrow \kk[z]$ given by $\sigma^{m}-1$. Let $G$ be the group of automorphisms generated by the following maps:
\begin{align*}
&\Theta_\beta: x \mapsto \beta x, \quad y \mapsto \beta\inv y, \quad z \mapsto z, \\
&\Psi_{m,\lambda}: 
	x \mapsto x, \quad 
    	y \mapsto y + \sum_{i=1}^n \frac{\lambda^i}{i!} \Delta_m^i(a)x^{im-1}, \quad 
    	z \mapsto z-m\lambda x^m, \\
&\Phi_{m,\lambda}: 
	x \mapsto x + \sum_{i=1}^n \frac{(-\lambda)^i}{i!} y^{im-1}\Delta_m^i(a), \quad 
	y \mapsto y, \quad 
    	z \mapsto z+m\lambda y^m.
\end{align*}
The automorphism group of $\kk[z](\sigma,a)$ is $G$ unless $a$ is \emph{reflective}, that is, there exists some $\rho \in \kk$ such that $a(\rho-z)=(-1)^n a(z)$  \cite[Theorem 3.29]{BJ}. When $a$ is reflective, the automorphism group contains an additional generator given by
\[ \Omega(x) = y, \quad \Omega(y) = (-1)^nx, \quad \Omega(z) = 1+\rho-z.\]
This work led naturally into considering the isomorphism problem for classical GWAs. In particular, $\kk[z](\sigma,a) \iso \kk[z](\sigma,a')$ if and only if $a'(z) = \eta a(\tau \pm z)$ for some $\eta,\tau \in \kk$ with $\eta \neq 0$ \cite[Theorem 3.28]{BJ}. This reflects an earlier remark\footnote{This was in a footnote, which you've totally been reading, right?} that $U_\lambda \iso U_{-\lambda}$.

The automorphism group of the first quantum Weyl algebra $\qwa$ was computed by Alev and Dumas \cite{AD1}. An adaptation of this result led to the corresponding isomorphism result: $\qwa \iso \pwa$ if and only if $p=q^{\pm 1}$ \cite{Giso1}.

As noted above, quantum Weyl algebras are examples of quantum GWAs. Let $R(\sigma,a)$ be a quantum GWA with $\sigma(z)=qz$, $q \in \kk^\times$. For the case $R=\kk[z^{\pm 1}]$ and $q$ is a nonroot of unity, Bavula and Jordan computed the automorphism group and solved the isomorphism problem \cite{BJ}. The case $R=\kk[z]$ and $q$ a nonroot of unity was completed by Richard and Solotar \cite{RS2}. The case for both base rings with $q$ a root of unity was completed by Su\'{a}rez-Alvarez and Vivas \cite{SAV}, see also \cite{KL3}. See \cite{CGWZ1} for a different approach, along with some results for higher rank quantum GWAs.

For higher rank quantized Weyl algebras, Goodearl and Hartwig determined the automorphism group and solved the isomorphism problem for generic parameters\footnote{Equivalently, the center of the algebra is $\kk$.} \cite{GH}.
This used a result of Jordan \cite{Jsimp}, which classified all height one prime ideals in the generic case\footnote{I told you this would pop up again.}. Levitt and Yakimov made use of \emph{noncommutative discriminants} to prove a corresponding result in the case that the center is a polynomial ring\footnote{Equivalent conditions on the parameters are given in the reference. It is useful to also mention that the algebra is free over its center in this case.} \cite{LY}. The case where the center is not trivial or a polynomial ring is open. See also the work of Hartwig and Rosso on graded isomorphisms of TGWAs over polynomial rings \cite{HR1}.

\subsection{Invariants}

We have already seen how in some ways the properties of the first Weyl algebra mimic those of the polynomial ring in two variables ($\GKdim\fwa = 2$) and differ ($\gldim\fwa=1$). Another crucial place where the two differ is with regards to the ring of invariants. If $G$ is a finite linear subgroup of $\Aut(\kk[x,y])$ acting without pseudo reflections, then $\kk[x,y]^G$ is another polynomial ring\footnote{In fact, this is an ``if and only if" statement due to the Shephard-Todd-Chevalley Theorem \cite{ShTo,Chev}.} and isomorphic again to $\kk[x,y]$. Remarkably\footnote{How remarkable depends on ones perspective.}, an invariant ring of $\fwa$ is \emph{never} isomorphic to $\fwa$. This was proved by Smith for solvable groups \cite{Sm1}. Alev, Hodges, and Velez proved the result for any finite group \cite{AHV}. Tikaradze showed that the first Weyl algebra is not the fixed ring of \emph{any} domain \cite{TIK}.

If the invariant ring of $\fwa$ is not another Weyl algebra, it is reasonable to ask what structure, if any, might be preserved. Suppose $\phi$ is an automorphism of a rank one GWA $R(\sigma,a)$ of order $\ell<\infty$ that restricts to an order $m$ automorphism of $R$, and there exists $\mu \in \kk^\times$ such that $\phi(x)=\mu\inv x$ and $\phi(y)=\mu y$. If $\gcd(\ell,m)=1$, then $(R(\sigma,a))^{\grp{\phi}} = R^{\grp{\phi}}(\sigma^\ell,\widehat{a})$ where 
\[ \widehat{a} = y^\ell x^\ell = \prod_{k=0}^{\ell-1} \sigma^{-k}(a).\]
In the case of $\fwa$, the invariant ring is a GWA. This result is due to Jordan and Wells in the case\footnote{In this case, the gcd condition is trivial.} that $m=1$ \cite[Theorem 2.6]{JW}. For the general statement, see \cite[Theorem 3.3]{GHo}. See \cite{GW1} for a discussion of invariants in the case of \emph{linear} automorphisms.

One can ask similar questions for other GWA-like structures. In doing this, one looks for a suitable set of conditions so that the fixed ring preserves the given structure. For example, for a BR algebra $B=R(t,\sigma,H,J)$ over a commutative noetherian ring, let $\phi \in \Aut(R)$ such that $\phi(H) \subset H$ and $\phi(J) \subset J$. For $\gamma \in \kk^\times$ there is an automorphism $\Phi_\lambda$ extending $\phi$ such that $\Phi_\gamma(t)=\gamma t$. If the order of $\gamma$ is $m<\infty$ and the order of $\phi$ is $n < \infty$, and if $\gcd(m,n)=1$, then the fixed ring by $\Phi_\lambda$ is another BR algebra \cite[Proposition 2.13]{GRW1}.

Now let $A=A_\mu(R,\sigma,a)$ be a TGWA of rank $n$ over a commutative domain $R$. Suppose that $\phi \in \Aut(A)$ restricts to an automorphism of $R$ of order $\ell < \infty$. Suppose that $\phi(X_i^{\pm})=\alpha_i^{\pm 1}$ for some $\alpha_i \in \kk^\times$ of order $m_i < \infty$. If $\ell$ and the $m_i$ are relatively prime, then $R^{\phi\mid_R}$ and the $(X_i^{\pm})^{m_i}$ generate another TGWA \cite[Theorem 3.10]{GR}. In certain cases, such as when $A$ is of type $(A_1)^n$, this TGWA is exactly the fixed ring $A^{\grp{\phi}}$. As a corollary this gives a result for higher rank GWAs.

\section{Future work and research directions}

We end this article by mentioning some ideas for future work in this area.

\subsection{Search for simple rings}
The Weyl algebra and its generalization give a host of examples of infinite dimensional simple rings. The diagram below shows the relationship between the various generalizations, where arrows represent inclusion\footnote{We have not included the class of quantizations as they do not fit nicely into the given paradigm.}.  To find further examples of infinite dimensional simple rings, it would be prudent to identify a suitable merger of higher rank BR algebras\footnote{We have not formally discussed higher rank BR algebras. A definition and discussion can be found in \cite{GRW2}.} and TGWAs.

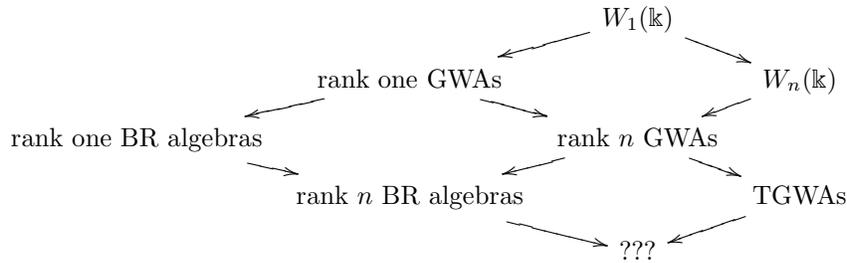
\vspace{-1em}
\[ \xymatrix@=0.4em{
& & \fwa \ar[dl] \ar[dr] & \\
& \text{rank one GWAs} \ar[dr] \ar[dl] & & \wan \ar[dl]  \\
 \text{rank one BR algebras} \ar[dr] & & \text{rank $n$ GWAs} \ar[dr]  \ar[dl] & \\
& \text{rank $n$ BR algebras} \ar[dr]   & & \text{TGWAs} \ar[dl] \\
&  & ??? &
}
\]

In \cite{GR3}, the authors consider \emph{twisted tensor products} of (T)GWAs and  those of BR algebras. This begins to give a picture of where the theories of TGWAs and BR algebras may intersect.

Not included in this diagram is the notion of a (rank one) \emph{weak generalized Weyl algebra} \cite{LMZweak}. In this setting, the map $\sigma$ is only required to be \emph{endomorphism}. A natural extension of this study would be to higher rank weak GWAs, as well as weak TGWAs and weak BR algebras.

\subsection{The Calabi-Yau condition}

Global dimension is connected with the Calabi-Yau (CY) condition \cite{ginz,RR2}, or more generally the twisted CY condition. We do not include the definition here, but note that $A_n(\kk)$ is CY. However, most proofs of this such as Berger's \cite[Theorem 6.5]{bergerCY}, relies on the fact that the associated graded ring of $A_n(\kk)$ is a polynomial ring\footnote{More precisely, it is the fact that $A_n(\kk)$ is a \emph{PBW deformation} of the polynomial ring in $2n$ variables.}, which itself is CY. Liu has proved that rank one GWAs over $\kk[z]$ and $\kk[z_1,z_2]$ are twisted CY \cite{liu,liu2}. It would be interesting to know whether GWAs over polynomial rings are twisted CY in general, as well as some more general criteria for this condition for GWAs, TGWAs, and BR algebras. See also the related study of Hochschild homology and cohomology for GWAs \cite{FSSA,SSAV}.

\subsection{A general isomorphism result}

As mentioned above, the isomorphism problem for higher rank quantized Weyl algebras is an open problem. In general, the \emph{graded isomorphism problem} is easier to solve. See, for example, the work of Hartwig and Rosso \cite{HR1}. Bell and Zhang proved that an isomorphism between connected $\NN$-graded algebras finitely generated in degree one implies an isomorphism as graded algebras \cite{BZ1}. A careful analysis of isomorphism problems of $\ZZ$-graded algebras \cite{BJ,SAV} suggests that a version of Bell and Zhang's isomorphism lemma should hold in the $\ZZ$-graded setting. However, the problem becomes a bit more subtle since there are isomorphisms that reverse the grading.

\subsection{Hopf actions}
The Weyl algebra and, more generally, GWAs exhibit \emph{quantum rigidity}, that is, their automorphism group is significantly smaller than their commutative counterparts. This becomes even more apparent when we restrict our attention to automorphisms that preserve the $\ZZ$-grading. Hence, for further (quantum) symmetries, it is natural to look for actions by finite-dimensional Hopf algebras. We refer to \cite{montgomery} for the basics of Hopf algebras and their actions on rings.

Cuadra, Etingof, and Walton proved that any action of a finite dimensional Hopf algebra $H$ on $A_n(\kk)$ factors through a group action \cite{CEW1,CEW2}. That is, if $H$ acts on $A_n(\kk)$, then there is a Hopf ideal $I$ of $H$ that acts trivially on $A_n(\kk)$ and $H/I$ is (isomorphic to) a group algebra. See \cite{EW4,EW3} for further results in this direction.

On the other hand, quantum GWAs $\kk[z](\sigma,a)$ with $\sigma$ finite order do exhibit quantum symmetry. In \cite{GW2}, it is shown that there are inner-faithful actions by \emph{generalized Taft algebras}. Moreover, these actions preserve the $\ZZ$-grading on these GWAs. See \cite{CWWZ} for results on filtered actions on certain quantum GWAs. The generalized Taft algebras are examples of pointed Hopf algebras. It would be interesting to know whether there are any actions by other families of Hopf algebras on quantum GWAs (other than group algebras).

\bibliography{gwabib}{}

\end{document}